\documentclass{article}
\usepackage[utf8]{inputenc}

\usepackage{amsthm}
\usepackage{amsfonts} 
\usepackage{amsmath}
\usepackage{mathtools}
\usepackage{geometry}
\usepackage{graphicx,amsmath,fullpage,mathptmx,helvet}
\usepackage[ruled, lined, linesnumbered, commentsnumbered, longend]{algorithm2e}
\usepackage{stmaryrd}
\usepackage{tikz}
\usepackage{tikz-qtree}
\usepackage{listings}
\usepackage{caption}
\usepackage{subcaption}
\usepackage{chngcntr}
\usepackage{authblk}
\usepackage{enumerate}
\usepackage{xurl}
\usepackage{mathrsfs}
\usepackage{pgfplots}

\usepackage[hidelinks]{hyperref}
\usepackage{todonotes}
\usepackage{cleveref}

\usetikzlibrary{shapes.geometric}

\theoremstyle{definition}
\newtheorem{definition}{Definition}[section]
\theoremstyle{plain}
\newtheorem{theorem}{Theorem}[section]
\newtheorem{example}{Example}[section]

\newtheorem{property}[theorem]{Property}

\newcommand*\diff{\mathop{}\!\mathrm{d}}

\DeclareMathOperator*{\argmax}{arg\,max}
\DeclareMathOperator*{\argmin}{arg\,min}
\DeclareMathOperator{\sublim}{sublim}
\theoremstyle{remark}
\newtheorem*{remark}{Remark}

\makeatletter
\@ifdefinable\skpt{\def\skpt#1\skpt{}}{}
\makeatother

\pgfplotsset{compat=1.18}

\title{Computational bounds on randomized algorithms for online bin stretching}
\author{Antoine Lhomme, Nicolas Catusse, Nadia Brauner}
\date{Univ. Grenoble Alpes, CNRS, Grenoble INP, G-SCOP, 38000 Grenoble, France\\[2ex]%
\today}

\begin{document}
\maketitle

\begin{abstract}
    A frequently studied performance measure in online optimization is \textit{competitive analysis}. It corresponds to the worst-case ratio, over all possible inputs of an algorithm, between the performance of the algorithm and the optimal offline performance. However, this analysis may be too pessimistic to give valuable insight on a problem. Several workarounds exist, such as \textit{randomized algorithms}.

    This paper aims to propose computational methods to construct randomized algorithms and to bound their performance on the classical online bin stretching problem. A game theory method is adapted to construct lower bounds on the performance of randomized online algorithms via linear programming. Another computational method is then proposed to construct randomized algorithms which perform better than the best deterministic algorithms known. Finally, another lower bound method for a restricted class of randomized algorithm for this problem is proposed.
\end{abstract}

\section{Online bin stretching}
The online bin stretching problem, introduced in \cite{AZAR200117}, is defined as follows: a finite sequence of items, each characterized by its size (represented by a real number), must be placed into $m$ bins. The items are revealed sequentially, meaning that an item, when revealed, must be placed into a bin before the next one is shown. The whole set of items given is also known to fit inside $m$ bins of unit size. 
The objective is to place all the items into the bins so that, in the end, the load of the fullest bin is minimized. An algorithm that places items into bins has its performance measured through its worst-case input; the performance of an algorithm is called in this problem the \textit{stretching factor}\footnote{If the sequence of items fits exactly into $m$ bins of unit size, then the stretching factor may also be seen as a competitive ratio.}.

This problem may also be interpreted as an online scheduling problem on $m$ parallel identical machines where the objective is to minimize the makespan (the total time to execute all tasks): using scheduling notations, the problem may be written as $P_m|online, OPT=1|C_{max}$, with $OPT=1$ indicating that the optimal makespan is known in advance to be of unit duration.

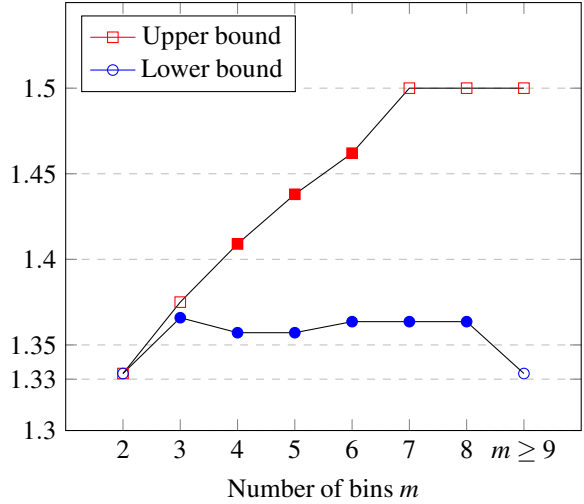
\begin{figure}[htbp]
    \centering
    \begin{tikzpicture}
    \begin{axis}[
        xlabel={Number of bins $m$},
        xmin=1, xmax=10,
        ymin=1.3, ymax=1.55,
        xtick={2, 3, 4, 5, 6, 7, 8, 9},
        xticklabels={2, 3, 4, 5, 6, 7, 8, $m\geq 9$},
        ytick={1.3, 1.33, 1.35, 1.4, 1.45, 1.5},
        legend pos=north west,
        ymajorgrids=true,
        grid style=dashed,
    ]
    
    \addplot[
        scatter,
        point meta=explicit symbolic, scatter/classes={
        a={mark=square,red}, b={mark=square*,red}
        }]
        coordinates {
        (2,1.3333) [a]
        (3,1.375) [a]
        (4,1.409) [b]
        (5,1.438) [b]
        (6,1.462) [b]
        (7,1.5) [a]
        (8,1.5) [a]
        (9,1.5) [a]
        };
    \addplot[
        scatter,
        point meta=explicit symbolic, scatter/classes={
        c={mark=o,blue}, d={mark=*,blue}
        }]
        coordinates {
        (2,1.3333) [c]
        (3,1.3659) [d]
        (4,1.35714) [d]
        (5,1.35714) [d]
        (6,1.3636) [d]
        (7,1.3636) [d]
        (8,1.3636) [d]
        (9,1.3333) [c]
        };
        \legend{Upper bound,,Lower bound,}
    \end{axis}
    
    \end{tikzpicture}
    \caption{Best known bounds on the minimal stretching factor of deterministic algorithms for the online bin stretching problem. Solid squares and circles represent bounds that were found through computational searches.}
    \label{obs_results}
\end{figure}

Interestingly, some of the best known lower and upper bounds (algorithms) for this problem in the deterministic\footnote{We employ the term \textit{deterministic} in contrast to \textit{randomized}: for some input, the output of a deterministic algorithm is  fixed, while the output of a randomized algorithm is a random variable} setting were recently found through computational searches: see \cite{Gabay2017}, \cite{BOHM20221}, \cite{lhomme2022online} for lower bounds and \cite{Liesk} for upper bounds. Figure~\ref{obs_results} depicts the current state of the art on this problem. Solid squares or solid circles correspond to values that were found through a computational search. Finding optimal online algorithms is still a challenge when the number of bins $m$ is greater than 2. 

\begin{figure}[htpb]
    \centering
    \resizebox{\linewidth}{!}{
	\begin{tikzpicture}
	    \tikzset{every tree node/.style={align=center,anchor=north}}
	    \tikzset{level distance=40pt}
	    \Tree
	    [.{$(0,0)$\\{\color{red} Next item : $1/3$}}
	    [.{$(1/3, 0)$\\{\color{red} Next item : $1/3$}}
	    [.{$(1/3, 1/3)$\\{\color{red} Next item : $1$}}
	    [.{$({\color{blue} 4/3}, 1/3)$} ]
	    [.{$(1/3, {\color{blue} 4/3})$} ]
	    ]
	    [.{$(2/3, 0)$\\{\color{red} Next item : $2/3$}}
	    [.{$({\color{blue} 4/3}, 0)$} ]
	    [.{$(2/3, 2/3)$\\{\color{red} Next item : $2/3$}}
	    [.{$({\color{blue} 4/3}, 2/3)$} ]
	    [.{$(2/3, {\color{blue} 4/3})$} ]
	    ]
	    ]
	    ]
	    [.{$(0, 1/3)$\\{\color{red} Next item : $1/3$}}
	    [.{$(1/3, 1/3)$\\{\color{red} Next item : $1$}}
	    [.{$({\color{blue} 4/3}, 1/3)$} ]
	    [.{$(1/3, {\color{blue} 4/3})$} ]
	    ]
	    [.{$(0, 2/3)$\\{\color{red} Next item : $2/3$}}
	    [.{$(0, {\color{blue} 4/3})$} ]
	    [.{$(2/3, 2/3)$\\{\color{red} Next item : $2/3$}}
	    [.{$({\color{blue} 4/3}, 2/3)$} ]
	    [.{$(2/3, {\color{blue} 4/3})$} ]
	    ]
	    ]
	    ]
	    ]
	\end{tikzpicture}
	}
    \caption{Adversarial strategy showing the lower bound of $\frac{4}{3}$ for the online bin stretching problem where 2 bins are available. In parentheses: occupied volume in each bin.}
    \label{tree_proof_4/3}
\end{figure}
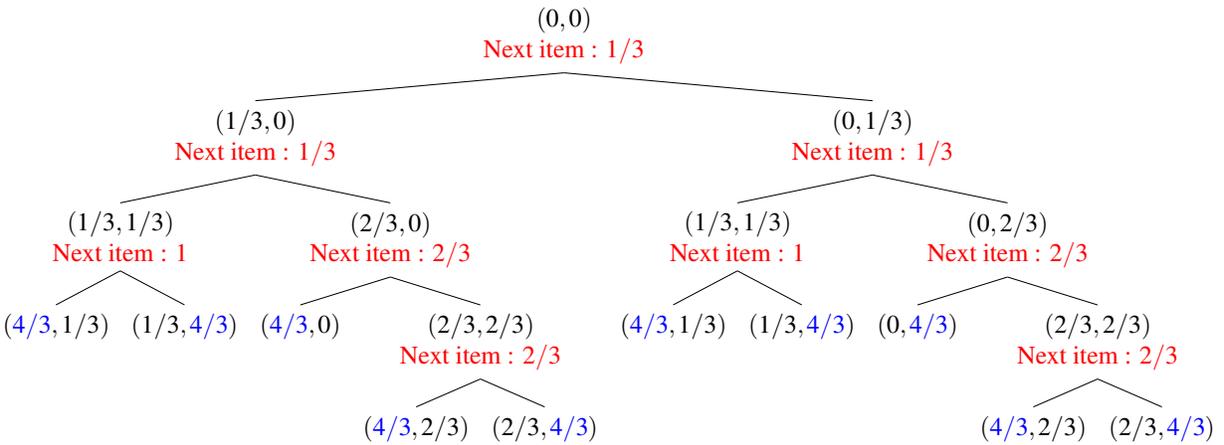

The computational method to find lower bounds may be summarized as follows: the problem is modelled as a two-player zero-sum game, where one player, the \textit{adversary}, sends items and the other player, the \textit{algorithm}, who places the items into bins. The algorithm aims to minimize the load of the largest bin. By restricting the adversary to a finite subset of possible item sizes, the game becomes finite - and hence its min-max value may be computed through standard techniques. Since the adversary player is restricted in this game, the min-max value is a lower bound of the best achievable competitive ratio of online deterministic algorithms. Figure~\ref{tree_proof_4/3} gives an example of an adversarial strategy that proves the lower bound $\frac{4}{3}$ when 2 bins are available.

To construct upper bounds, the same kind of technique is used. The problem is once again modelled as a game between the adversary and the algorithm; however, to ensure that a strategy for the algorithm player corresponds to a true online bin stretching algorithm, the adversary player sends item \textit{classes} instead of items, which correspond to an interval of possible item sizes. See \cite{Liesk} for more details about this technique, which is not explained here.

The heart of the computational methods to obtain the lower and upper bounds consists of min-max searches in a game tree. Naturally, the size of the game tree becomes exponentially greater as one allows more items to be sent, and computation times are the limiting factor for obtaining better bounds.

Most of the results known so far for this problem are in the deterministic setting. The aim of this paper is to explore the randomized setting and to do so with computational methods, in the same fashion as to what has been done for the deterministic case. \Cref{sec:Randomized_algorithms} justifies the interest of randomized algorithms and describes a framework to study online problems. In \Cref{sec:lower_bounds_with_pl}, a linear programming methods allows the computational search of lower bounds for online problems. This method is used with success in \Cref{LB on random alg for obs} on the online bin stretching problem to construct new lower bounds on the performance of randomized algorithms; moreover \Cref{proof convergence} details a proof of the convergence of that method for the online bin stretching problem. Moreover, a restricted class of randomized algorithms is studied in \Cref{sec:finding_randomized_algorithms} and \Cref{sec:lower_bounds_on_M2}, which yields a randomized algorithm that outperforms deterministic algorithms in the simple case for two bins.

\section{Randomized algorithms}\label{sec:Randomized_algorithms}
\subsection{Motivation}\label{sec:Motivation}
Consider the following game of \textit{matching pennies}: two players, A and B, each have a coin. Each picks a side, heads or tails, and then both show each other their coin. If both players have selected the same side, then player A wins. Otherwise, player B wins.
A deterministic strategy for a player corresponds to choosing heads or tails, thus there are only two possible deterministic strategies. Trivially, the worst-case scenario for any given deterministic strategy of any player is a loss, which makes the worst-case analysis not so insightful. A well known work-around is to consider \textit{randomized strategies}, which are informally, strategies making use of randomness in some way. For games with a quantified objective, the worst-case performance of a randomized strategy can be defined as the worst-case (relative to the strategies of the other players) of the expected pay-off (relative to the randomness of the randomized strategy) of the game.

In the previous example of matching pennies, the randomized strategy that picks heads (otherwise tails) with probability $\frac{1}{2}$ is the best one in the sense of worst-case analysis: regardless of the strategy chosen by the other player, the likelihood of winning is $\frac{1}{2}$.

In online optimization, considering randomized algorithms is a fairly standard practice (see \cite{Ben-David1994}) as doing so helps to lessen the flaws of worst-case analysis, as showcased in the previous matching pennies example. Thus, we aim to consider randomized algorithms for the previously introduced online bin stretching problem. So far, the only result known on randomized algorithms for online bin stretching is the following lower bound, from \cite{Gabay2017}:
\begin{property}
    For the online bin stretching problem with two bins, no randomized algorithms may have a stretching factor lower than $7/6$.
\end{property}

In game theory, two types of randomized strategies that are often considered are \textit{mixed strategies} and \textit{behavioral strategies}. The first one corresponds to picking, before the start of the game, a deterministic strategy according to some probability distribution and then following the chosen deterministic strategy for the rest of the game. The second one corresponds to, whenever a player has to make a move, picking a move according to some distribution. It has been shown in \cite{kuhn} that, if the game is finite and of \textit{perfect recall}, \textit{i.e.},  when players "remember" their past moves, then both types of strategies are equivalent. If the game is infinite, then it has also been shown in \cite{Aumann+1964+627+650} that both kinds of strategies are equivalent, however this requires more complex definitions of mixed and behavioral strategies: these will be addressed in \Cref{sec:definition_randomized_algorithms}.

\subsection{The request-answer framework}\label{section RA framework}

Let us first remark that the whole idea of using randomized algorithms requires \textit{independence} between the decisions taken by an algorithm and the input sequence received. For instance, the input may be the current weather and the decision may be whether one purchases an umbrella or not. It is then reasonable to assume that the decision taken does not have any influence on the weather for future upcoming days. Our formalism is built on this assumption: it corresponds to \textit{request answer games} \cite{Borodin_el_yaniv}, which we will briefly define in the following.

As the performance of an algorithm is measured in this paper by its worst-case input, an \textit{adversary} player will be in charge of sending the input to the algorithm with the goal of making the algorithm perform as poorly as possible - this forms a two-player, zero-sum game.

\begin{definition}[Request answer game]
    A request answer game $\mathcal G = (\mathcal I, \mathcal D, f)$ consists of:
    \begin{itemize}
        \item A set of instances $\mathcal I$, where instances are finite sequences of \textit{requests};
        \item A set of decisions $\mathcal D$;
        \item A pay-off function $f\colon F \to \mathbb R$, where $F = \{(I, (d_1, \dots, d_k))~|~I\in\mathcal I~\text{ and } (d_1, \dots, d_k)\in\mathcal D^{|I|}\}$.
    \end{itemize}
\end{definition}

The game proceeds as follows: a request arrives to the algorithm, who has to pick some decision in $\mathcal D$. This is repeated, but the sequence of requests is constrained to be in the set of instances $\mathcal I$. When the sequence ends, the pay-off is computed through the pay-off function $f$ - in all of this document, the algorithm always aims to minimize the pay-off. For the online bin stretching problem, an instance is a sequence of items that fit into $m$ bins of unit size, a decision is the index of a bin into which to pack an item, and the pay-off function corresponds to the load of the fullest bin.

The objective in the deterministic setting is to find the best algorithm possible in the worst-case; formally, we aim to compute the \textit{deterministic value} $v^*_{det}$ of the game:
$$v_{det}^*=\inf_{A \in \mathcal A_{det}}\;\sup_{I\in\mathcal I} f(I, A(I))$$
In the above equation, $\mathcal A_{det}$ represents the set of deterministic algorithms, and we suppose that $A(I)$ corresponds to the sequence of decisions taken by the algorithm on the input sequence $I$.

In the randomized case, the output $A(I)$ of an algorithm $A$ on some instance $I$ is a random variable; thus we aim to compute the \textit{randomized value} of the game defined as:
$$v_{rand}^*= \inf_{A \in \mathcal A_{rand}}\;\sup_{I\in\mathcal I} \mathbb E(f(I, A(I)))$$

\begin{remark}
    Online algorithms are often studied through \textit{competitive analysis}, which measures the \textit{competitive ratio} of algorithms: it corresponds to the worst-case over all inputs of the ratio of performance between the pay-off from the online algorithm and an optimal \textit{offline} algorithm, \textit{i.e.}, who has knowledge of the entire input in advance. With the given formalism, one may define the pay-off function to correspond to such a ratio, so that the value of the game corresponds to the best competitive ratio achievable by an algorithm.
\end{remark}

\section{Lower bounds for online randomized algorithms with linear programming}\label{sec:lower_bounds_with_pl}
In this section, we use a game theory technique to compute lower bounds on the performance of randomized online algorithms for a given online problem. The method is inspired from \cite{Stengel} and uses linear programming. The following develops both the general idea alternatively with a practical toy example, which is then applied in Section~\ref{LB on random alg for obs} to the online bin stretching problem to construct new bounds. 

The method we present here is inspired by the work of \cite{Stengel}, in which a method to compute the equilibrium of a game via a linear program is given. It is adapted here for online problems; the method cannot be applied straightforwardly to the game tree of the request-answer game, as the whole input sequence must be chosen in advance before the algorithm makes random choices.

Let $(\mathcal I, \mathcal D, f)$ be a request-answer game as defined in Section~\ref{section RA framework}, with $\mathcal D$ finite and $\mathcal I$ finite as well. We aim to compute the worst case expected pay-off of the best randomized algorithm for this request-answer game. We require $\mathcal I$ to be finite here, which is very restrictive. 
In practice, one may approximate $\mathcal I$ by finite subsets $(\mathcal I_g)_{g\in\mathbb N}$, such that $\forall g\in\mathbb N,~\mathcal I_g \subseteq \mathcal I_{g+1} \subset \mathcal I$. Then, applying the method to the game $(\mathcal I_g, \mathcal D, f)$ will result in a lower bound on the performance of randomized algorithms, since the adversary is restricted to a subset of instances. Increasing $g$ may then give better and better bounds.

\begin{example}
To clarify the idea, the following finite two-player zero-sum game, represented with its game tree in Figure~\ref{Game example figure}, will serve as an example. The first player, named the \textit{algorithm}, aims to go skiing during the weekend. Each day, the algorithm may either go skiing (\textbf{S}) or stay at home, drinking hot chocolate (\textbf{C}). Each day, it will either be snowing (\textbf{S}) or the weather will be clear (\textbf{C}). The algorithm aims to minimize some unhappiness function defined as follows: skiing while it is snowing for a day gives 1 unhappiness, but -2 unhappiness in clear weather. Moreover, if it was snowing on the first day, going skiing in clear weather yields -3 unhappiness thanks to fresh that day. Drinking hot chocolate while the weather is clear gives 1 unhappiness, but it gives -1 unhappiness while it is snowing. The weather is, however, completely unpredictable in this particular ski resort and is only known when the algorithm has already made its choice. The worst-case analysis of this problem is not very insightful: if the algorithm decides to go skiing, then it may end up snowing, and if the algorithm stays at home, the weather may be clear; the min-max value of the game tree in Figure~\ref{Game example figure} is 2. In other words, the best deterministic algorithm has a worst-case of 2. 

\begin{figure}[htpb]
    \centering
    \begin{tikzpicture}
        [level distance=10mm,
        level 1/.style={sibling distance=60mm},
        level 2/.style={sibling distance=30mm},
        level 3/.style={sibling distance=15mm},
        level 4/.style={sibling distance=7mm}]
        \node [draw, diamond]{}
        child  {
            node [fill, circle] {}
            child  {
                node [draw, diamond] {}
                child  {
                    node [fill, circle] {}
                    child  {
                        node [draw, rectangle] {2}
                        edge from parent node[left] {S}
                    }
                    child  {
                        node [draw, rectangle] {-2}
                        edge from parent node[right] {C}
                    }
                    edge from parent node[left] {S}
                }
                child  {
                    node [fill, circle] {}
                    child  {
                        node [draw, rectangle] {0}
                        edge from parent node[left] {S}
                    }
                    child  {
                        node [draw, rectangle] {2}
                        edge from parent node[right] {C}
                    }
                    edge from parent node[right] {C}
                }
                edge from parent node[left] {S}
            }
            child  {
                node [draw, diamond] {}
                child  {
                    node [fill, circle] {}
                    child  {
                        node [draw, rectangle] {-1}
                        edge from parent node[left] {S}
                    }
                    child  {
                        node [draw, rectangle] {-4}
                        edge from parent node[right] {C}
                    }
                    edge from parent node[left] {S}
                }
                child  {
                    node [fill, circle] {}
                    child  {
                        node [draw, rectangle] {-3}
                        edge from parent node[left] {S}
                    }
                    child  {
                        node [draw, rectangle] {-1}
                        edge from parent node[right] {C}
                    }
                    edge from parent node[right] {C}
                }
                edge from parent node[right] {C}
            }
          edge from parent
            node[above left] {S}
        }
        child  {
            node [fill, circle] {}
            child  {
                node [draw, diamond] {}
                child  {
                    node [fill, circle] {}
                    child  {
                        node [draw, rectangle] {0}
                        edge from parent node[left] {S}
                    }
                    child  {
                        node [draw, rectangle] {-4}
                        edge from parent node[right] {C}
                    }
                    edge from parent node[left] {S}
                }
                child  {
                    node [fill, circle] {}
                    child  {
                        node [draw, rectangle] {-2}
                        edge from parent node[left] {S}
                    }
                    child  {
                        node [draw, rectangle] {0}
                        edge from parent node[right] {C}
                    }
                    edge from parent node[right] {C}
                }
                edge from parent node[left] {S}
            }
            child  {
                node [draw, diamond] {}
                child  {
                    node [fill, circle] {}
                    child  {
                        node [draw, rectangle] {2}
                        edge from parent node[left] {S}
                    }
                    child  {
                        node [draw, rectangle] {-1}
                        edge from parent node[right] {C}
                    }
                    edge from parent node[left] {S}
                }
                child  {
                    node [fill, circle] {}
                    child  {
                        node [draw, rectangle] {0}
                        edge from parent node[left] {S}
                    }
                    child  {
                        node [draw, rectangle] {2}
                        edge from parent node[right] {C}
                    }
                    edge from parent node[right] {C}
                }
                edge from parent node[right] {C}
            }
          edge from parent
            node[above right] {C}
        }
        ;
    \end{tikzpicture}
    \caption{\label{Game example figure}Two player zero-sum game example. Diamond nodes represent a decision node for the algorithm player (aiming to minimize the pay-off) and filled circles are for the adversary player (aiming to maximize the pay-off). S represents "Snow" for the adversary and "Skiing" for the algorithm, while C represents "Clear weather" for the adversary and "Chocolate" for the algorithm.}
\end{figure}
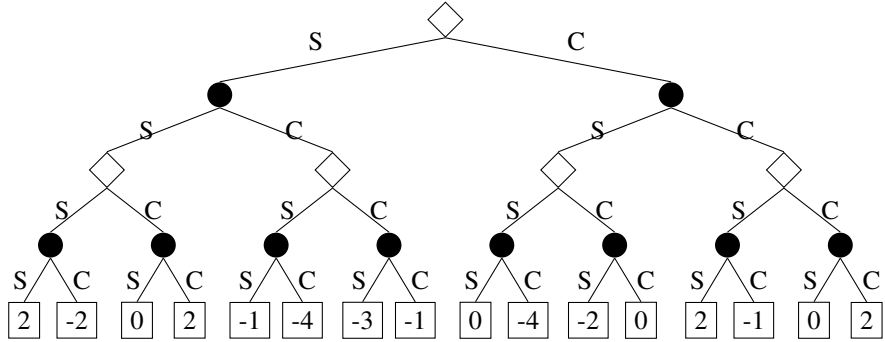

\end{example}

To compute the expected value of the game for the optimal randomized algorithm in the worst case, we model the game as a linear program in the following way.

First, label all possible moves of the algorithm in this game as $x_i$, and label all the sequences of moves of the adversary as $y_j$. These will be variables in our linear program that represent the randomized strategies of the algorithm and adversary player. For the algorithm, variables $x_i$ will represent a \textit{flow} of probability - this may be clearer in the example that follows.

Flow constraints are added as follows: if $x_{f_1}, \dots, x_{f_k}$ are variables corresponding to the algorithm moves for some game state, and that $x_p$ is the variable corresponding to the last move made by the algorithm, then add the constraint: 
$$\sum_{i=1}^k x_{f_i} = x_p$$
For the very first move of the algorithm player, replace $x_p$ by 1. 
These constraints in matrix form may be written as:
$$Ex = e$$

For the adversary, add a variable $y_i$ for each possible sequence of moves for the adversary, \textit{i.e.}, a variable for each element of $\mathcal I$. This variable represents the probability for the adversary of choosing such a move sequence. Then, add the constraint:
$$\sum_{i}y_i = 1$$

Moreover, $x$ and $y$ are constrained to be positive.

\begin{example} 
Figure~\ref{Game example figure labelled} shows a possible labelling of the edges for the game of Figure~\ref{Game example figure}. 
From the root node, the algorithm has to make a decision. The probability of choosing the left branch is exactly $x_1$, and is $x_6$ for the right branch. From this node's left side child, the algorithm has to make its second move. It picks left with probability $\frac{x_2}{x_1}$ and right with probability $\frac{x_3}{x_1}$.

The adversary has four sequences of moves available: (\textbf{S}, \textbf{S}), (\textbf{S}, \textbf{C}), (\textbf{C}, \textbf{S}) and (\textbf{C}, \textbf{C}).
Thus, variables $y_1, y_2, y_3$ and $y_4$ are added. 

The flow constraints for the algorithm player can be written as:
\begin{align*}
    &x_1 + x_6 = 1\\
    &x_2+x_3 = x_1\\
    &x_4+x_5 = x_1\\
    &x_7+x_8 = x_6\\
    &x_9+x_{10} = x_6
\end{align*}
The constraint for the adversary player is written as:
\begin{align*}
    &y_1 + y_2 + y_3 + y_4 = 1
\end{align*}
These constraints may then be written in the form $Ex = e$ and $\sum_i y_i = 1$, with $E$ and $e$ defined as:
$$E = \begin{pmatrix}
    1 & 0 & 0 & 0 & 0 & 1 & 0 & 0 & 0 & 0\\
    -1 & 1 & 1 & 0 & 0 & 0 & 0 & 0 & 0 & 0\\
    -1 & 0 & 0 & 1 & 1 & 0 & 0 & 0 & 0 & 0\\
    0 & 0 & 0 & 0 & 0 & -1 & 1 & 1 & 0 & 0\\
    0 & 0 & 0 & 0 & 0 & -1 & 0 & 0 & 1 & 1\\
\end{pmatrix},~~~e = \begin{pmatrix}
    1\\
    0\\
    0\\
    0\\
    0\\
\end{pmatrix}  $$

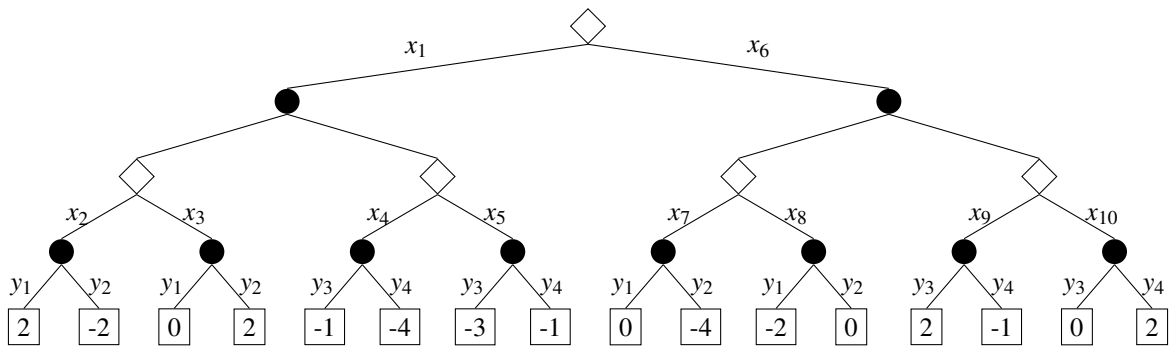
\begin{figure}[htpb]
    \centering
    \begin{tikzpicture}
        [level distance=10mm,
        level 1/.style={sibling distance=80mm},
        level 2/.style={sibling distance=40mm},
        level 3/.style={sibling distance=20mm},
        level 4/.style={sibling distance=10mm}]
        \node [draw, diamond]{}
        child  {
            node [fill, circle] {}
            child  {
                node [draw, diamond] {}
                child  {
                    node [fill, circle] {}
                    child  {
                        node [draw, rectangle] {2}
                        edge from parent node[left] {$y_1$}
                    }
                    child  {
                        node [draw, rectangle] {-2}
                        edge from parent node[right] {$y_2$}
                    }
                    edge from parent node[left] {$x_2$}
                }
                child  {
                    node [fill, circle] {}
                    child  {
                        node [draw, rectangle] {0}
                        edge from parent node[left] {$y_1$}
                    }
                    child  {
                        node [draw, rectangle] {2}
                        edge from parent node[right] {$y_2$}
                    }
                    edge from parent node[right] {$x_3$}
                }
                edge from parent node[left] {}
            }
            child  {
                node [draw, diamond] {}
                child  {
                    node [fill, circle] {}
                    child  {
                        node [draw, rectangle] {-1}
                        edge from parent node[left] {$y_3$}
                    }
                    child  {
                        node [draw, rectangle] {-4}
                        edge from parent node[right] {$y_4$}
                    }
                    edge from parent node[left] {$x_4$}
                }
                child  {
                    node [fill, circle] {}
                    child  {
                        node [draw, rectangle] {-3}
                        edge from parent node[left] {$y_3$}
                    }
                    child  {
                        node [draw, rectangle] {-1}
                        edge from parent node[right] {$y_{4}$}
                    }
                    edge from parent node[right] {$x_5$}
                }
                edge from parent node[right] {}
            }
          edge from parent
            node[above left] {$x_1$}
        }
        child  {
            node [fill, circle] {}
            child  {
                node [draw, diamond] {}
                child  {
                    node [fill, circle] {}
                    child  {
                        node [draw, rectangle] {0}
                        edge from parent node[left] {$y_{1}$}
                    }
                    child  {
                        node [draw, rectangle] {-4}
                        edge from parent node[right] {$y_{2}$}
                    }
                    edge from parent node[left] {$x_7$}
                }
                child  {
                    node [fill, circle] {}
                    child  {
                        node [draw, rectangle] {-2}
                        edge from parent node[left] {$y_{1}$}
                    }
                    child  {
                        node [draw, rectangle] {0}
                        edge from parent node[right] {$y_{2}$}
                    }
                    edge from parent node[right] {$x_8$}
                }
                edge from parent node[left] {}
            }
            child  {
                node [draw, diamond] {}
                child  {
                    node [fill, circle] {}
                    child  {
                        node [draw, rectangle] {2}
                        edge from parent node[left] {$y_{3}$}
                    }
                    child  {
                        node [draw, rectangle] {-1}
                        edge from parent node[right] {$y_{4}$}
                    }
                    edge from parent node[left] {$x_9$}
                }
                child  {
                    node [fill, circle] {}
                    child  {
                        node [draw, rectangle] {0}
                        edge from parent node[left] {$y_{3}$}
                    }
                    child  {
                        node [draw, rectangle] {2}
                        edge from parent node[right] {$y_{4}$}
                    }
                    edge from parent node[right] {$x_{10}$}
                }
                edge from parent node[right] {}
            }
          edge from parent
            node[above right] {$x_6$}
        }
        ;
    \end{tikzpicture}
    \caption{\label{Game example figure labelled}Two player zero-sum game example. Diamond nodes represent a decision node for the algorithm player and filled circles are for the adversary player. Variables $x$ and $y$ have been added to represent probabilities for possible moves.}
\end{figure}

\end{example}

In our context of worst-case analysis, observe that the adversary has no interest in playing a randomized strategy: for any strategy for the player algorithm, there exists a deterministic strategy for the adversary which is optimal (see \cite{kuhn}); however writing these equations in such a way allows the constraints to be linear rather than integer.

Remark that the probability of reaching some leaf is exactly the product between the last $x_j$ variable and the $y_i$ variable corresponding to the adversary move sequence. In the previous example, see Figure~\ref{Game example figure labelled}, the probability of reaching the leftmost leaf is exactly $x_2\cdot y_1$. Given a strategy $x$ and a strategy $y$ for each player, the expected value may be computed through the weighted sum of the pay-off of each leaf where weights are the probability of reaching each leaf. Hence, there exists a matrix $C$ such that for any $x$ and $y$ strategies for both players, the expected value of the game can be written as $x^T Cy$.

\begin{example}
In the previous example, the matrix $C$ is written as:
$$C = \begin{pmatrix}
    0 & 0 & 0 & 0\\
    2 & -2 & 0 & 0\\
    0 & 2 & 0 & 0\\
    0 & 0 & -1 & -4\\
    0 & 0 & -3 & -1\\
    0 & 0 & 0 & 0\\
    0 & -4 & 0 & 0\\
    -2 & 0 & 0 & 0\\
    0 & 0 & 2 & -1\\
    0 & 0 & 0 & 2
\end{pmatrix}$$
And the expected value of the game, given two randomized strategies $x$ and $y$ can hence be written as $x^T C y$.

\end{example}
The following optimization problem yields the best worst-case expected pay-off of the game:

$$\min_x \max_y x^T C y$$
such that:
$$Ex = e,\;\;\sum_i y_i=1,\;\;x\geq 0,\;\;y\geq 0$$

Let $x$ be any strategy for the algorithm. Then, the best response for the adversary can be obtained by solving the problem:

$$\max_y x^TCy$$
such that:
$$\sum_i y_i=1,\;\;y\geq 0$$

As in \cite{Stengel}, we dualize this problem into:

$$\min_u u$$
such that:
$$u \begin{pmatrix} 1\\ \vdots \\ 1 \end{pmatrix}  \geq C^Tx$$

By strong duality, it has the same optimal value as its dual. So the previous non-linear problem may be written as the following linear program:

$$\min_{u, x} u$$
such that:
$$Ex = e,\;\;u \begin{pmatrix} 1\\ \vdots \\ 1 \end{pmatrix}  \geq C^Tx,\;\;x\geq 0$$

Solving this linear program yields the best randomized algorithm $x$ and its worst-case expected value $u$, which corresponds to the best worst-case expected pay-off of a randomized algorithm for the game.

Generating this linear program can be done by exploring the game tree through an exhaustive search.

\begin{example}
    
With the previously given matrices, the linear program may be solved with any linear optimization solver. The optimal randomized worst-case found is $u \approx -0.48$, with $x_7$ and $x_9$ being 0 in an optimal solution. The best deterministic algorithm had a worst-case pay-off of 2.

\end{example}

One may observe that the game tree should always begin with an algorithm decision and end with an adversary decision. If the game tree starts with an adversary decision, then the adversary always pick the best subtree for itself. Similarly, if the last move is an algorithm move, then the algorithm always pick the best move for itself.

\section{Bounds on randomized algorithms for online bin stretching}

\subsection{Lower bounds with the linear programming method}\label{LB on random alg for obs}
The linear programming technique exposed in Section~\ref{sec:lower_bounds_with_pl} was applied successfully to the online bin stretching problem. This allowed the construction of the following bounds:

\begin{theorem}
    For $m=2$ bins, there is a lower bound on the competitive ratio of randomized algorithm of $\frac{7}{6} \approx 1.1667$. For $m=3$, there is a lower bound of $\frac{19}{16} = 1.1875$, and for $m=4$, of $\frac{11}{9} \approx 1.2222$.
\end{theorem}

\begin{remark}
    The lower bound of $\frac{7}{6}$ was already known from \cite{these_gabay}, but was found through a different kind of analysis.
\end{remark}

As explained in Section~\ref{sec:lower_bounds_with_pl}, the linear programming method requires the game to be finite. To do so, we restrict the adversary to items of the form $x/g$ with $x\in \{1, \dots, g\}$ with $g$ some integer fixed in advance. A move of the adversary corresponds to choosing one item to send, and a move of the algorithm corresponds to choosing a bin into which the given item should be placed. Since the adversary is restricted in its choice of items, finding the optimal algorithm of that game with the linear program yields a lower bound on the performance of randomized algorithms.

Constructing the linear program may be done in a straightforward manner, \textit{i.e.}, by exploring the whole game tree. The implementation of our code may be found online\footnote{\url{https://gricad-gitlab.univ-grenoble-alpes.fr/lhommea/computational-proofs-on-randomized-algorithms-for-online-bin-stretching}}, which also contains a visual representation of a game tree for online bin stretching alongside the associated weights.

Our linear program depends on some parameter $g\in\mathbb N$, such that the adversary is restricted to items of the form $\frac{x}{g}$ with $x\in\{1, \dots, g\}$. As such, increasing $g$ makes the adversary player stronger, and may result in better bounds at the cost of having a bigger game tree. A natural question is then: do the lower bounds parametrized by $g$ converge towards the optimal performance of a randomized online algorithm? It may be the case that they do not - we give in Section~\ref{sec:theoretical_limit} examples of games where such a convergence is not possible in a broader setting. However, the convergence holds for the online bin stretching problem:

\begin{theorem}
    The method to construct lower bounds on randomized algorithms does converge for the online bin stretching problem.
\end{theorem}

The proof is available in Appendix~\ref{proof convergence}.

\subsection{Finding randomized algorithms through computational searches}\label{sec:finding_randomized_algorithms}

So far, our method managed to bound the performance of randomized algorithms - however, as yet, no randomized algorithms that outperform deterministic algorithms for the online bin stretching problem are known. Is there a randomized algorithm for the online bin stretching problem that performs strictly better than any deterministic one? 

To answer this, we consider a subset of randomized algorithms which consists of \textit{mixed} algorithms that are comprised of only two deterministic algorithms. This class will be henceforth written as $\mathcal M_2$. An algorithm of this class, before receiving an input sequence, selects one of its two deterministic algorithms according to some probability and then takes the same decisions as the chosen deterministic algorithm for the whole sequence.

In \cite{Liesk}, a computational method was proposed to construct online deterministic algorithms through a min-max search in a game tree; we propose an adaptation of that method to construct randomized algorithms from the class $\mathcal M_2$. When only two bins are available, the best competitive ratio for deterministic algorithms is exactly $\frac{4}{3}$; our aim is to find a pair of deterministic algorithms $A_1$ and $A_2$ together with a probability $p$ so that regardless of the instance $I$, $p(f(I, A_1(I))) + (1-p)(f(I, A_1(I))) <  \frac{4}{3}$.

\begin{figure}[htpb]
    \centering
    \begin{tikzpicture}
        [level distance=12mm,
        level 1/.style={sibling distance=52mm},
        level 2/.style={sibling distance=52mm},
        level 3/.style={sibling distance=32mm},
        level 4/.style={sibling distance=15mm}]
        \node {$p=\frac{1}{2}$}
            child {node [draw, rectangle] {$((4, 1), (3, 2))$}
                child {
                node [draw, circle] {  }
                    child {
                        node [draw, rectangle] {\begin{tabular}{c} $((5, 1), (4, 2))$ \\ Pay-off : $9/2$ \end{tabular}}
                    }
                    child {
                        node [draw, rectangle] {\begin{tabular}{c} $((5, 1), (3, 3))$ \\ Pay-off : $4$ \end{tabular}}
                    }
                    child {
                        node [draw, rectangle] {\begin{tabular}{c} $((4, 2), (4, 2))$ \\ Pay-off : $4$ \end{tabular}}
                    }
                    child {
                        node [draw, rectangle] {\begin{tabular}{c} $((4, 2), (3, 3))$ \\ Pay-off : $7/2$ \end{tabular}}
                    }
                edge from parent node[above left] {Next item class: $1$}
                }
                child [dashed] {}
            }
        ;
    \end{tikzpicture}
    \caption{\label{Randomized Upper bound method illustration}Illustration of the game tree where two algorithms play in parallel. The pay-off at a leaf is the weighted sum according to $p$ of the highest bin loads of each algorithm.}
\end{figure}
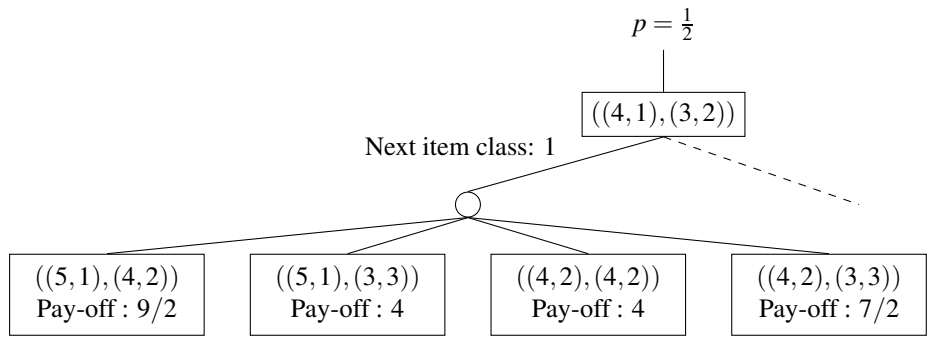

\cite{Liesk} models the problem of finding deterministic algorithms for the online bin stretching problem as a finite two player game; in that game, the \textit{adversary} player sends items and the algorithm places the items into bins, but the items sizes are only known within some interval. The model depends once again on a discretization setting; the higher it is, the smaller the item intervals are. Finding a strategy for the algorithm player in this game does yield a valid online bin stretching algorithm. The game tree is illustrated in Figure~\ref{Upper bound method illustration} - for the sake of clarity, the details of that method are omitted here: the actual game tree is slightly more complex\footnote{Since item classes and bin loads represent intervals, adding an item into some bin may result in a wider interval - the idea is to keep the interval size constant and to branch on the possible resulting intervals. The game tree is actually wider than what is shown in Figure~\ref{Upper bound method illustration} and Figure~\ref{Randomized Upper bound method illustration}. The technique is explained in details in \cite{Liesk}.}.

We propose to modify that game so that the algorithm plays two moves (not necessarily distinct) in parallel, which corresponds to having two deterministic algorithms at once (illustrated in Figure~\ref{Randomized Upper bound method illustration}). Each game state holds the bin loads of both algorithms. Then, the pay-off function at a leaf node is the weighted sum of the pay-off of both algorithms at that leaf node, according to some probability $p$ fixed beforehand of choosing either algorithm.

\begin{figure}[htpb]
    
    \centering
    \begin{tikzpicture}
        [level distance=10mm,
        level 1/.style={sibling distance=52mm},
        level 2/.style={sibling distance=26mm},
        level 3/.style={sibling distance=13mm},
        level 4/.style={sibling distance=7mm}]
        \node {}
            child {node [draw, rectangle] {$(4, 1)$}
                child {
                node [draw, circle] {  }
                    child {
                        node [draw, rectangle] {$(5, 1)$}
                        child[dashed] { node {}}
                        child[dashed] { node {}}
                    }
                    child {
                        node [draw, rectangle]{$(4, 2)$}
                        child[dashed] { node {}}
                        child[dashed] { node {}}
                    }
                edge from parent node[above left] {Next item class: $1$}
                }
                child {
                node [draw, circle] {  }
                    child {
                        node [draw, rectangle] {$(6, 1)$}
                        child[dashed] { node {}}
                        child[dashed] { node {}}
                    }
                    child {
                        node [draw, rectangle]{$(4, 3)$}
                        child[dashed] { node {}}
                        child[dashed] { node {}}
                    }
                edge from parent node[above right] {Next item class: $2$}
                }
            }
        ;
    \end{tikzpicture}
    \caption{\label{Upper bound method illustration}Illustration of the game tree for constructing upper bounds in the deterministic case. In parentheses, the load of each bin. Circle nodes correspond to game states where the algorithm has to take a decision.}
\end{figure}
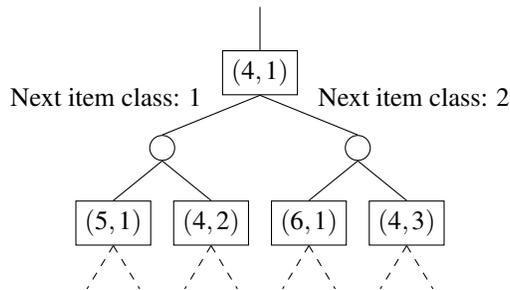

This idea, however, is fairly inefficient as we are essentially iterating over all possible pairs of deterministic algorithms - and the probability $p$ needs to be chosen. We propose two improvements to speed up the method and reach higher discretization settings in order to obtain better algorithms:

\begin{enumerate}
    \item When $p = \frac{1}{2}$, both algorithms may not be distinguished. To reduce the number of possible game states, one may hence always order the loads of the bins of the first algorithm to be greater lexicographically than the loads of the bins of the second algorithm. 
    \item In \cite{lhomme2022online}, a propagation technique was presented that allows a strong pruning of the game states for the deterministic, lower bound case. The idea is to consider partial information (only bin loads) of a game state to quickly check if the algorithm player is certain to win from that game state. This may be adapted straightforwardly to the upper bound game, and to our method where two algorithms play at once, by considering the $2m$-uplet of bin loads as the partial information. 
\end{enumerate}

Moreover, classical optimization techniques for min-max searches may be used in conjunction to these, such as hash tables to keep visited nodes in memory and a dynamic program to quickly compute which items may be sent at any configuration (which is a NP-complete problem, see \cite{lhomme2022online} for more technical details on the dynamic program).

The implementation of that search leads to the following result:
\begin{theorem}
    There is a $\frac{5}{4}$-competitive randomized algorithm for the online bin stretching problem with $2$ bins, which is better than the optimal deterministic algorithm (which is $\frac{4}{3}$-competitive).
\end{theorem}

This result was obtained on with discretization 18, \textit{i.e.}, when considering only item sizes of the form $\frac{x}{18}$ with $x$ integer, with roughly 35 seconds of computation time on a machine with Intel Core i5-1135G7 and 15 GB of RAM.

\subsection{Lower bounds on the performance of the $\mathcal M_2$ class}\label{sec:lower_bounds_on_M2}

So far, we have a method to construct lower bounds and a method for upper bounds, but a significant gap still exists even for the seemingly simple case $m=2$: the best algorithm in the deterministic case has competitive ratio $4/3\approx 1.334$, the best randomized known has competitive ratio $1.25$ but the best lower bound known on randomized algorithm is only $7/6\approx 1.1667$.

The same idea of having two algorithm players play in parallel may be applied to do the same for lower bounds; in other words, we aim to compute lower bounds on the performance of randomized algorithms from the class $\mathcal M_2$ that are distributions over 2 deterministic algorithms. 

An idea could be to fix some probability $p\in(0,1)$ corresponding to the probability of selecting the first deterministic algorithm of the pair, and then do a search while all algorithms play simultaneously, with the score of a leaf being defined as the weighted sum of the score of each algorithm. However, doing so results in a lower bound proof only for \textbf{one single} probability $p\in(0, 1)$. To bridge this gap, one may observe that if the pay-off function of the online problem is bounded, then a slightly weaker proof is valid for a neighborhood of $p$. For online bin stretching, the pay-off function is bounded since it corresponds to the load of the fullest bin; as the item sequence is guaranteed to fit inside $m$ bins of unit size, the load of the fullest bin is at most $m$ (putting all items into the same bin). This is formalized through the following property.

\begin{property}
    Let $p\in(0,1)$ and $v\geq 0$ such that:
    $$\forall (A_1, A_2)\in\mathcal (A_{det})^2,~ \forall I\in \mathcal I,\;\;~ p\cdot f(I, A_1(I)) + (1-p)\cdot f(I, A_2(I)) \geq v$$
    Then:
    $$\forall p'\in (0,1),~\forall (A_1, A_2)\in\mathcal (A_{det})^2,~ \forall I\in \mathcal I,~ \;\;p'\cdot f(I, A_1(I)) + (1-p')\cdot f(I, A_2(I)) \geq (v - |p-p'|\cdot m)$$ 
\end{property}

\begin{proof}
    The term $p\cdot f(I, A_1(I)) + (1-p)\cdot f(I, A_2(I))$ is linear in $p$ and its derivative is at most $|f|\leq m$.
\end{proof}

In other words, when finding a lower bound $v$ for some probability $p$, one also has lower bounds of values close to $v$ for a probability in the neighborhood of $p$. This property leads to the following scheme to derive lower bounds on the performance of randomized algorithms that are a probability distribution over two deterministic algorithms, which corresponds to finding lower bounds for several probabilities $p$ and then using the above property to have a general lower bound on \textbf{any} probability. 

Suppose that we have a \textit{minmax} search function which takes as input some probability $p$ and some target lower bound to be proved, and that returns true if the lower bound was proved on any algorithm from $\mathcal M_2$ with probability $p$. Then Algorithm~\ref{alg:cap} presents how one may construct a lower bound on the performance of any algorithm from $\mathcal M_2$. Note that we only need to consider $p$ within $(0, \frac{1}{2})$, as if $p$ is greater than $\frac{1}{2}$ then $1-p$ is smaller than $\frac{1}{2}$ and plays the same role.

\begin{algorithm}[htpb]
    \caption{A scheme to construct lower bounds on the performance of algorithms from the class $\mathcal M_2$}\label{alg:cap}
    \KwIn{Discretization $g\in\mathbb N$, number of min-max searches $N\in\mathbb N$, $t\geq 0$ the lower bound to be proved}
    \KwOut{True if the given lower bound is proved, False otherwise}
    $\delta \gets \frac{1}{4N}$ \tcp*[f]{The distance between two tested probabilities will be 2$\delta$}

    \For{$i = 0, \dots, N-1$} {
        $p \gets \delta + 2*\delta * i$

        \If{minmax$(g, t + m\delta, p)$ == False} { 

        \tcp*[f]{The min-max function is called here with a slightly stronger bound, $t + m\delta$,than the one we aim to prove, and returns true if that bound was proved when the probability of any algorithm from $\mathcal M_2$ is fixed to $p$}

            \KwRet{False}
        }
        
    }
    \KwRet{True}
\end{algorithm}

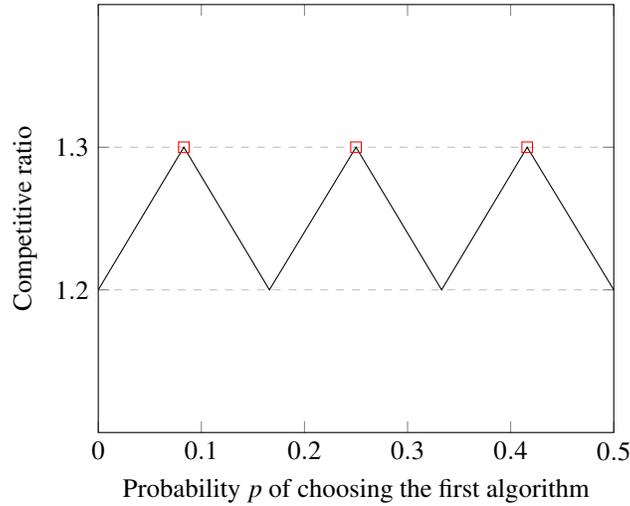
\begin{figure}[htbp]
    \centering
    \begin{tikzpicture}
    \begin{axis}[
        xlabel={Probability $p$ of choosing the first algorithm},
        ylabel={Competitive ratio},
        xmin=0, xmax=0.5,
        ymin=1.1, ymax=1.4,
        xtick={0, 0.1, 0.2, 0.3, 0.4, 0.5},
        ytick={1.2, 1.3},
        legend pos=south east,
        ymajorgrids=true,
        grid style=dashed,
    ]
    
    \addplot[
        scatter,
        point meta=explicit symbolic, scatter/classes={
        c={mark=square,red}, d={mark=circle,blue}
        }]
        coordinates {
        (0, 1.2)[d]
        (0.083,1.3)[c]
        (0.166,1.2)[d]
        (0.25,1.3)[c]
        (0.333,1.2)[d]
        (0.416,1.3)[c]
        (0.5,1.2)[d]
        };
    \end{axis}
    
    \end{tikzpicture}
    \caption{Illustration of the methods to construct lower bounds on the class $\mathcal M_2$. Finding the lower bounds for a fixed $p$ (squares) implies weaker bounds on a neighborhood.}
    \label{Figure algo M2}
\end{figure}

The application of this scheme to our problem yields the following result:
\begin{theorem}
    No randomized algorithm from $\mathcal M_2$ has a competitive ratio smaller than $29/24\approx 1,208$ for the online bin stretching when 2 bins are available.
\end{theorem}
This result directly bounds the performance of the upper bound method to construct randomized algorithm, which yielded an algorithm with competitive ratio $1.25$ when two bins are available. 

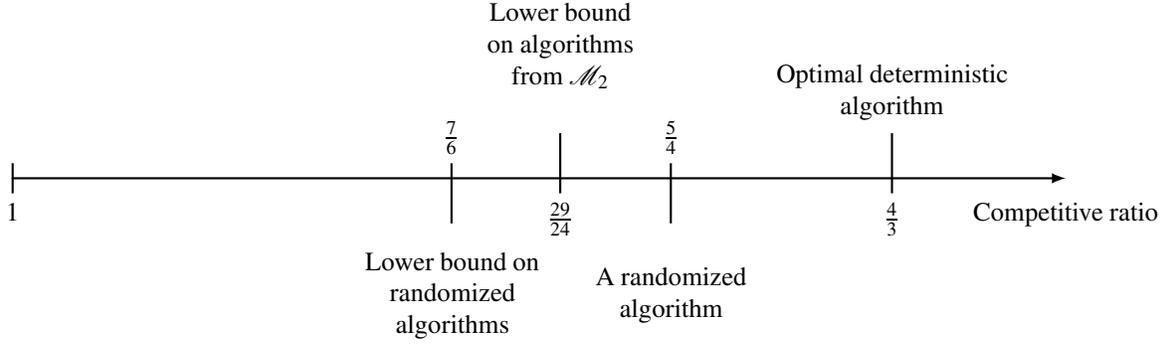
\begin{figure}
\begin{center}

\begin{tikzpicture}[
    scale = 10,
    x=3.5cm,
    legend/.style={
        below = 0.2cm
    },
]
    \draw[black,->,thick,>=latex,line cap=rect]
        (0,0) -- (0.4,0);

    \draw (0,0) node [legend] {$1$};
    \draw[black,thick]
        (0,-0.02) -- (0,0.02);

    \draw (0.1667,0) node[above=0.2cm] {$\frac{7}{6}$}; 
    \draw[black,thick]
        (0.1667,-0.06) -- (0.1667,0.02);

    \draw (0.208,0) node[legend] {$\frac{29}{24}$};
    \draw[black,thick]
        (0.208,-0.02) -- (0.208,0.06);

    \draw (0.25,0) node[above = 0.2cm] {$\frac{5}{4}$};
    \draw[black,thick]
        (0.25,-0.06) -- (0.25,0.02);

    \draw (0.334,0) node[legend] {$\frac{4}{3}$};
    \draw[black,thick]
        (0.334,-0.02) -- (0.334,0.06);

    \draw (0.4,0) node[legend] {Competitive ratio};

    \draw (0.334,0) node[above = 0.6cm] {\begin{tabular}{c}Optimal deterministic\\algorithm\end{tabular}};

    \draw (0.25,0) node[below = 1cm] {\begin{tabular}{c} A randomized\\ algorithm\end{tabular}};

    \draw (0.208,0) node[above = 1cm] {\begin{tabular}{c}Lower bound\\on algorithms\\from $\mathcal M_2$\end{tabular}};

    \draw (0.167,0) node[below = 0.8cm] {\begin{tabular}{c} Lower bound on\\randomized\\algorithms\end{tabular}};

\end{tikzpicture}
\caption{\label{fig:bounds_found}Bounds for online bin stretching when $m$=2 bins are available.}
\end{center}
\end{figure}

\Cref{fig:bounds_found} presents our results for online bin stretching when two bins are available. Once again, our code is available online here\footnote{ \url{https://gricad-gitlab.univ-grenoble-alpes.fr/lhommea/computational-proofs-on-randomized-algorithms-for-online-bin-stretching}}.

\subsection{A theoretical limit of the computational approach}\label{sec:theoretical_limit}

One may define, for any $k\in\mathbb N$, the class $\mathcal M_k$ of mixed algorithms consisting of $k$ deterministic algorithms. However, doing so may not necessarily improve algorithms: the following gives two zero-sum two-player game, where one player, the \textit{algorithm}, aims to minimize the pay-off and the other player, the \textit{adversary} aims to maximize the pay-off. In these games, the algorithm player has a deterministic worst-case of 1 and a randomized worst-case expected pay-off of 0; however if the algorithm uses a mixed strategy that is a distribution on a countable set of deterministic strategies, then the worst-case expected pay-off becomes 1. This shows that it may not be possible to approach a 'good' randomized algorithm through a distribution on a growing finite number of deterministic algorithms.

\textbf{First game example:}
 Two players each pick a positive real number, respectively denoted $x$ and $y$ for the algorithm and the adversary. The pay-off function of the game is then defined as: $$f(x,y) = \begin{cases} \frac{1}{\sqrt{|x-y|}} \text{ if }x\neq y\\0 \text{ otherwise} \end{cases}$$
 The algorithm aims to minimize $f$ while the adversary tries to maximize it. The game may be interpreted as both players picking a number, with the adversary trying to guess the number chosen by the algorithm as closely as possible. That game does not have a bounded worst-case pay-off, as the adversary may pick $y$ close to $x$ to reach an arbitrarily high value, so any deterministic strategy for the algorithm has a worst-case of $+\infty$. Similarly, if the algorithm strategy is randomized over a countable set of deterministic strategies, the adversary may choose $y$ close to some possible output $x$ of the algorithm randomized strategy to reach an arbitrarily high expected pay-off, so the worst-case expected pay-off of any randomized strategy defined as a probability over a countable set of deterministic strategies has a worst-case expected pay-off of $+\infty$. However, the randomized strategy $A$ that follows the distribution $p(x)= e^{-x}$ has bounded worst-case expected pay-off:
 $$\forall y\in\mathbb R,\;\mathbb{E}(f(A, y)) = \int_{x\in\mathbb R} p(x) f(x, y)\diff x \leq  2\sqrt\pi$$
 
 This example, however, makes use of an unbounded pay-off function. It is hence natural to ask if a bounded pay-off function is sufficient to ensure that any randomized strategy may be approached by mixed strategies using countable deterministic strategies. However, the following example shows that this is not the case.

 \textbf{Second game example:}
Both players play simultaneously then both moves are revealed, and the pay-off is given by a pay-off function $f$ defined below. 
Let $D = \{0, 1, \dots, 9\}$ be the set of digits. A move of any player is picking some infinite sequence of digits. One may see a move as a real number between 0 and 1. Let $M = D^\mathbb N$ denote the set of moves. The adversary aims, once again, to guess the move of the algorithm. However, the adversary is here able to guess a countable set of moves: if the algorithm chosen number is within that sequence, then the adversary wins. This is formalized as follows: 

Define $b\colon M^\mathbb N \to M$ to be some bijection between the set of sequences of moves and the set of moves. Such a bijection exists because $\mathbb R$ and $\mathbb R^\mathbb N$ have the same cardinality.

The pay-off resulting from move $x = (x_1, x_2, \dots) \in M$ by player 1 and move $y=(y_1, y_2, \dots)\in M$ of player 2 is defined as:
\[
    f(x,y)= 
\begin{cases}
    1& \text{if } \exists (z^1, z^2, \dots, )\in M ^\mathbb N,~\exists k\in\mathbb N \\ & \text{such that }b(z^1, z^2, \dots, z^k, x, z^{k+1}, z^{k+2}, \dots) = y\\
    0              & \text{otherwise}
\end{cases}
\]
The bijection $b$ allows the adversary to guess a countable set of moves. If the guess is successful, \textit{i.e.}, if the move from the algorithm is within the guessed sequence from the adversary, then the adversary wins (pay-off of 1), otherwise the algorithm wins (pay-off of 0).

If player 1 has a mixed strategy that is a distribution on a countable set of deterministic strategies $(x^1, x^2, \dots) \in M^\mathbb N$, then there is a move $y$ for player 2 such that the expected pay-off is 1: fix $y = b(x^1, x^2, \dots)$. Then:
$$\forall i \in \mathbb{N},~f(x^i,y) = 1$$
Thus, the expected pay-off is 1.

Let us show that player 1 has a randomized strategy that has an expected pay-off of 0 regardless of the strategy of player 2. If player 1 picks each digit in his sequence randomly, then:
$$\forall y\in M,~\mathbb E_x (f(x,y)) \leq \sum_{k\in \mathbb{N}} P_k$$
with:
$$P_k = \mathbb P(\exists (z^1, z^2, \dots, )\in M ^\mathbb N, ~|~b(z^1, z^2, \dots, z^k, x, z^{k+1}, z^{k+2}, \dots) = y)$$

However: 
$$P_k = 0$$

Hence, the expected pay-off is 0.

These examples show that, in a general setting, it may not be possible to approximate a randomized strategy by considering mixed strategies using countable deterministic strategies. This may be a limit to the power of computational methods.

\section{Conclusion and perspectives}
The competitive analysis framework to study online algorithms is often pointed out for its flaws. Workarounds exist, such as randomized algorithms; however this paradigm is more complex than deterministic algorithms - thus this paper proposes methods for a computational approach to perform competitive analysis on randomized algorithms. First, a method derived from game theory is proposed to compute the best randomized online algorithm in a finite setting via linear programming. This method was then used with success on the online bin stretching problem to construct lower bounds on the performance of randomized algorithms. Moreover, the adaptation of a computational method to construct upper bounds in the deterministic case yields a new randomized algorithm for online bin stretching that outperform the best known one in the simple case where 2 bins are available. Besides, we also propose an algorithm to give lower bounds on a simple class of randomized algorithms.

These results show that the paradigm of considering randomized algorithms may be studied through computational tools. Several questions arise next: is it possible to find an efficient pruning of the game tree for the linear programming method? Can a more efficient method be found to compute online randomized algorithms? Could computational methods be extended to other related online problems, such as online scheduling problems, or even to other online optimization paradigms, such as online algorithms with predictions?

\appendix

\section{On the definition of randomized online algorithms}\label{sec:definition_randomized_algorithms}
As stated in Section~\ref{sec:Motivation}, two types of strategies are prevalent in game theory: \textit{mixed} strategies and \textit{behavioral} strategies. A mixed strategy, before the start of the game, selects a deterministic strategy according to some probability distribution, and then applies the chosen deterministic strategy for the whole game. A behavioral strategy of a player contains the probability distribution of selecting a move for each possible game state where that player has to make a move. Both types of strategies have been shown to be equivalent in \cite{kuhn} when the game is finite and players have \textit{perfect recall}, \textit{i.e.}, remember their past moves.

However, if the game has an infinite number of moves and if players have a continuity of decisions available to them, the definition itself of a mixed strategy is not so clear. For example, given two deterministic strategies $A_1$ and $A_2$, the mixed strategy that picks $A_1$ with probability $\frac{1}{2}$ and $A_2$ otherwise is well-defined. But in an infinite game, there may be an infinite number of deterministic strategies; the set of deterministic strategies may not even be countable. To define a distribution on this space may not be easy, as this requires one to be able to measure the space of deterministic strategies. Some works such as \cite{Ben-David1994} or \cite{renault2014lower} consider randomized strategies to be a distribution over the space of deterministic strategies, which, as pointed out in \cite{Aumann+1964+627+650}, leads to several complications. Indeed, the space of deterministic strategies may be fairly complex; constructing a measure on that space may not even be possible. To circumvent this difficulty, \cite{Aumann+1964+627+650} proposes new definitions for mixed and behavioral strategies and then shows that both types of strategies are equivalent. Since both are equivalent, this paper only makes use of the behavioral strategy definition from \cite{Aumann+1964+627+650}:

\begin{definition}[Behavioral strategy]
    A behavioral strategy $b$ for some player is a sequence of measurable transformations $b_1, b_2, \dots$ such that $b_i \colon X_i \times \Omega_i \to Y_i$, where $X_i$ is the input space prior to the $i$-th move, $\Omega_i$ is a sample space and $Y_i$ is the set of decisions that may be taken. Given some input $x\in X_i$, the function $~b_i(x, \cdot)$ corresponds to the random variable associated with the action taken at move $i$ when the input is $x$. The functions $(b_i)_{i\in\mathbb N}$ must be measurable transformations, \textit{i.e.}, $\forall B_i\subseteq Y_i,~B_i\text{ measurable }$ then $b_i^{-1}(B_i)$ must be measurable. This condition is to ensure that if all players have a randomized strategy then the expected pay-off of the game may be computed. 
\end{definition}

\section{Proof of convergence of the linear program for online bin stretching}\label{proof convergence}
The lower bounds method for randomized algorithms in the online bin stretching case does depend on some discretization (or \textit{granularity}) parameter. 
This section aims to prove that the lower bounds obtained by the aforementioned method do converge towards the optimal competitive ratio of randomized algorithm.

\begin{theorem}
    The method to construct lower bounds on randomized algorithms does converge for the online bin stretching problem.
\end{theorem}

\textbf{Remark:} The statement of the theorem might be slightly misleading, in that we do not show that the sequence of lower bounds does converge towards the optimal value when the discretization increases, but rather we show that this sequence has a subsequence that converges towards the optimal value.

\textbf{Notations for the proof:} 

\begin{itemize}
    \item $\mathcal I$ represents the set of finite item sequences that fit into $m$ bins of unit size.
    \item $\mathcal I_g$ is the set of finite item sequences that fit into $m$ bins of unit size where all item sizes are of the form $\frac{x}{g}$ with $x$ integer.
    \item $\mathcal A$ is the set of randomized algorithms. We use the definition given in \Cref{sec:definition_randomized_algorithms} for a randomized algorithm. 
    \item $\mathcal D = \{d_1, \dots, d_m\}$ corresponds to the set of decisions of the algorithm: $d_i$ is the decision to place an item in the $i$-th bin. 
    \item $f$ is the pay-off function. It takes as input a finite sequence of items and a finite sequence of decisions with the same length, and outputs the load of the fullest bin after placing the items in the corresponding bins.
    \item $\mu\colon \mathcal P(\Omega)\to \mathbb R_{\geq 0}$ is the measure function on the sample space $\Omega$. We consider that all sample spaces are copies of the interval $[0, 1]$.
    \item A randomized online algorithm $A$ is a sequence of measurable functions $A_1, A_2, \dots$ such that the function $A_i\colon X_i \times \Omega_i \to \mathcal D$ represents the decision process for the $i$-th move of the algorithm; the set $X_i$ to the information space of the algorithm prior to taking the $i$-th decision. An element from this set corresponds to the input of the algorithm before taking the $i$-th decision, \textit{i.e.}, the sequence of objects received and decisions taken so far alongside a new item to be placed. An element from $X_i$ is hence written as $(I, D)$ where $I\in \mathcal I$, $|I| = i$ and $D\in \mathcal D^{i-1}$. For a given input $x = (I, D)\in X_i$, the function: \begin{align*}
        \Omega_i &\to \mathcal D \\
        \omega & \mapsto A_i(x, \omega)
      \end{align*}
      corresponds to the random variable of the decision taken. For a given input $x\in X_i$ and a decision $d\in \mathcal D$, the probability of taking the decision $d$ at step $i$ is hence:
      $$\mathbb P_{\omega\in\Omega_i} (A_i(x, \omega) = d) = \mu(A_i(x, \cdot)^{-1} (d))$$
\end{itemize}

\textbf{Proof idea:} 
Solving the linear program for some discretization $g$ yields a vector $x$ representing an algorithm strategy for a game where item sizes are restricted to fractions of the form $\frac{n}{g}$, with $n$ integer. The idea of the proof is to construct a randomized algorithm $A^{\infty, \epsilon}$ from the algorithms obtained with the linear program method, such that this constructed algorithm may process any item sequence from $\bigcup_{k\in\mathbb{N}} \mathcal I_{2^k}$  and does not perform worse (up to $\epsilon$) than the highest lower bound obtained through the method.

That algorithm is not yet a true online bin stretching algorithm, as it only deals with items of a form $\frac{x}{2^k}$ with $x$ and $k$ integers; however, we show that we may "extend" the algorithm into $A^{\infty, \epsilon, \delta}$ in order to process any item sequence from $\mathcal I$ without changing the competitive ratio up to an arbitrarily small constant $\delta$. In other words, we construct an algorithm that has performance arbitrarily close to the best lower bound found through our method; since the performance of any algorithm is also an upper bound on the optimal performance of an algorithm, this implies the convergence of the method.

Before starting the proof, the following definition is proposed:
\begin{definition}[Sublimit]
    Define the \textit{sublimit} of a sequence as the set of all limits of convergent subsequences of this sequence. Formally, if $(u_n)_{n \in \mathbb{N}}$ is a sequence, then:
    $$ \underset{n\to +\infty}{\sublim}\; u_n = \{ \lim_{n \to +\infty} v_n \text{ such that } (v_n)_{n \in \mathbb{N}} \text{ is a convergent subsequence of } (u_n)_{n \in \mathbb{N}}  \}$$
\end{definition}
If a sequence takes values in a compact set, then Bolzano-Weierstrass Theorem ensures that its sublimit is non-empty.

\begin{proof}
     Denote by $v^*_g$ the lower bound obtained by the method when the discretization is $g$ and denote by $A^g$ the corresponding randomized algorithm. Define $v^*$ to be the optimal competitive ratio of a randomized algorithm:
    $$v^* = \inf_{A\in\mathcal A} \sup_{I\in \mathcal I} \mathbb E (f(I, A(I)))$$
    
    We will show that: $$\lim_{k\to +\infty} v^*_{2^k} = v^*$$ 
    This proof only consider powers of 2 as discretization, to ensure that $\mathcal I_{2^k} \subseteq \mathcal I_{2^{k+1}}$, \textit{i.e.}, that the adversary has an increasing set of possible moves; this implies that the sequence $(v_{2^k}^*)$ is increasing (as the adversary is more and more powerful) and converging, since the sequence is bounded by the performance of any randomized algorithm.

    The proof is done in three steps: 
    \begin{enumerate}
        \item Construct $A^{\infty, \epsilon}$, an online algorithm that may only process item sequences from $\bigcup_{k\in\mathbb{N}} \mathcal I_{2^k}$
        \item Prove that $A^{\infty, \epsilon}$ has a competitive ratio of at most $\lim\limits_{k\to+\infty} v^*_{2^k}+\epsilon$
        \item Construct $A^{\infty, \epsilon, \delta}$ an online algorithm that has competitive ratio at most $v^*_{2^k}+\epsilon+\delta$
    \end{enumerate}

    \textbf{1. Construction of $A^{\infty, \epsilon}$}:

    Let $\epsilon>0$. We construct $A^{\infty, \epsilon} = (A^{\infty, \epsilon}_1, A^{\infty, \epsilon}_2, \dots)$ an online algorithm that may process any item sequence from $\bigcup_{k\in\mathbb{N}} \mathcal I_{2^k}$, such that the algorithm is $\lim\limits_{k\to +\infty} v^*_{2^k} +\epsilon$ competitive. By definition of behavioral algorithms given in \Cref{sec:definition_randomized_algorithms}, an online algorithm must be some sequence of measurable transformations: functions $A^{\infty, \epsilon}_i$ are defined on $X_i \times \Omega_i $ and map to $\mathcal D$, where $\Omega_i$ is a sample space, $X_i$ is the input space of the algorithm at step $i$, \textit{i.e.}, the item sequence given so far and the $i-1$ decisions taken so far. 

    Functions $A_i^{\infty, \epsilon}$, as well as sets $G(I)$ for any input sequence $I$ of length $i$ are constructed inductively on $i$. Define first $G(\emptyset) = \mathbb N$, with $\emptyset$ denoting here the empty item sequence, and define as well $\Delta_i = \frac{\epsilon}{m^2 2^i}$.

    Let $i\in \mathbb N^*$. We aim to construct $A^{\infty, \epsilon}_i$.
    First, let us observe that the set of sequences of decisions of length $i-1$ is finite; as such it is possible to give an ordering $D_1, \dots, D_{m^{i-1}}$ of $\mathcal D^{i-1}$.

    Let $p^j_{i, D, d}(I)$ be the probability for algorithm $A^j$ to select move $d\in\mathcal D$ when the sequence of items received is $I \in \bigcup_k \mathcal I_{2^k}$ (with $I$ of length $i$) and the sequence of decisions taken so far is $D\in \mathcal D^{i-1}$. Formally:
    $$p^j_{i, D, d}(I) = \mathbb P_{\omega\in \Omega_i}(A^j_i((I, D), \omega) = d)$$ 
    Remark that such a probability is not necessarily well-defined for any $j\in\mathbb N$: the algorithm $A^j$ may only deal with items of the form $\frac{x}{j}$ with $x$ integer.

    Define then $P^j_{i}(I)$ to be the $m^i$-uplet containing all $p^j_{i, k, d}(I)$ for $d\in\mathcal D$ and $k\in\{1, \dots, m^{i-1}\}$:

\begin{align*}
    P^j_{i}(I) = \Bigl( & p^j_{i, D_1, d_1}(I), p^j_{i, D_1, d_2}(I), \dots, p^j_{i, D_1, d_m}(I),\\
    & p^j_{i, D_2, d_1}(I), p^j_{i, D_2, d_2}(I), \dots, p^j_{i, D_2, d_m}(I),\\
    & \dots\\
    & p^j_{i, D_{m^{i-1}}, d_1}(I), p^j_{i, D_{m^{i-1}}, d_2}(I), \dots, p^j_{i, D_{m^{i-1}}, d_m}(I) \Bigl)
\end{align*}

    Then, one may observe that the $m^i$-uplet $P^j_i(I)$ is in the space $H_i$: 
    
    \begin{align*}
        P^j_i(I) \in H_i = \Bigl\{(x_{s})\in\mathbb R^{m^i}~| & ~\forall s\in\{1, \dots, m^{i}\}~x_s\in (0, 1)\\
        & \forall k\in\{0, \dots, m^{i-1}-1\},~\sum\limits_{l=1}^m x_{mk + l} = 1 \Bigl\}
    \end{align*}

    $H_i$ is compact; as such, Bolzano-Weierstrass Theorem ensures that any sequence from this set, and in particular the sequence $(P^j_i(I))_{j\in\mathbb N}$, has a subsequence that converges within $H_i$. The idea of our construction is to pick such a limit for the probabilities of the algorithm $A^{\infty, \epsilon}$. Since several limits may exist we need to have a procedure to select one. This may be done either immediately with the axiom of choice, or by remarking that the sublimit (the set of limits of convergent subsequences) is itself compact: since $H_i$ is compact, hence bounded, the sublimit of $(P^j_i(I))_{j\in\mathbb N}$ is clearly bounded; moreover, the sublimit is closed, since a converging sequence of limits is also a limit. A way to pick a singular element from a finite dimensional compact set is to order it with respect to some lexical ordering, and to select the minimum. The minimum always belongs to the compact.

    Define: $$P^{\infty, \epsilon}_i(I) = \min\limits_{\leq_{lex}} \underset{j\in G(I|_{1, \dots, i-1})}{\sublim} P^j_i(I)$$

    In the above equation, $\leq_{lex}$ denotes any fixed lexical ordering on $\mathbb R^{m^i}$, and $I|_{1, \dots, i-1}$ corresponds to the first $i-1$ items received.

    Define also $p^{\infty, \epsilon}_{i, D, d}$ the corresponding value from the $m^i$-uplet $P^{\infty, \epsilon}_i(I)$ for a given $D\in \mathcal D^{m^{i-1}}$ and $d\in\mathcal D$.

    Let:
    $$G(I) = \{j\in G(I|_{1, \dots, i-1})~|~||P^{\infty, \epsilon}_i(I) - P^{j}_i(I)||_{\infty} \leq \Delta_i\} $$
    The set $G(I)$ corresponds to indices of algorithms $A^g$ that are \textit{close} from the algorithm $A^{\infty, \epsilon}$ constructed so far; close means that the probability for taking some decision at step $i$ is distant by at most $\Delta_i$. By definition of the sublimit, $G(I)$ is always infinite. 

    Finally, the function $A^{\infty, \epsilon}_i$ may be defined as follows:

    $$A^{\infty, \epsilon}_i((I, D), \omega) = d_l \text{ where }l \text{ is such that }\omega\in \left] \sum_{s=1}^{l-1}p^{\infty, \epsilon}_{i, D, d_s}(I), \sum_{s=1}^{l}p^{\infty, \epsilon}_{i, D, d_s}(I) \right] $$

    In other words, the probability for the algorithm $A^{\infty, \epsilon}$ to select decision $d$ at move $i$ given input $(I, D)$ is exactly $p^{\infty, \epsilon}_{i, D, d}$.

    Remember that so far, we only defined $A^{\infty, \epsilon}$ on item sequences from $\bigcup_{k \in \mathbb N}I_{2^k}$ which is a countable set (since each $I_{2^k}$ is finite); as such we can show that $A^{\infty, \epsilon}_i$ is a measurable function:
    $$\forall d\in\mathcal D,~\;(A^{\infty, \epsilon}_i)^{-1}(d) = \underbrace{\bigcup_{\substack{I\in \bigcup_{k\in\mathbb N}\mathcal I_{2^k}\\ |I| = i}}\;\bigcup_{D\in\mathcal D^{i-1}} ((I, D), (\underbrace{A^{\infty, \epsilon}_i((I, D), \cdot))^{-1}(d)}_{\text{segment by construction}})}_{\text{Countable union}} $$

    \textbf{2. Proof that $A^{\infty, \epsilon}$ has competitive ratio at most $\lim_k v_{2^k}^* + \epsilon$}:

    Let $I\in\bigcup_{k\in\mathbb N}\mathcal I_{2^k}$. We will show that the expected score (load of the fullest bin) of the algorithm $A^{\infty, \epsilon}$ when receiving the item sequence $I$ is close from the expected score of the algorithm $A^j$ on $I$ for any $j\in G(I)$.

    By construction, given some input $x\in X_i$, the probability of choosing some move $d$ at step $i$ for the algorithm $A^{\infty, \epsilon}$ is close from the probability of choosing the same move with the same input $x$ of any algorithm from $G(I)$:
    $$\forall x\in X_i,~\forall d\in\mathcal D,~\;\; ||\mathbb P_{\omega\in\Omega_i}(A^{\infty, \epsilon}_i(x, \omega) = d) - \mathbb P_{\omega\in\Omega_i}(A^{j}_i(x, \omega) = d)||_{\infty} \leq \Delta_i$$

    Let $I\in\bigcup_{k\in\mathbb N}\mathcal I_{2^k}$ with $|I| = s$. Let $A^{\infty, \epsilon}(I)$ be the random variable corresponding to the output of the algorithm when receiving input $I$, which takes values in $\mathcal D^{s}$. The expected score of $A^{\infty, \epsilon}$ on $I$ may be written as:

    $$ \mathbb E(f(I, A^{\infty, \epsilon}(I))) = \sum_{D\in\mathcal D^s} \mathbb P(A^{\infty, \epsilon}(I) = D)f(I, D) $$
    $$ = \sum_{D\in\mathcal D^s} \prod_{i=1}^s \underbrace{p^{\infty, \epsilon}_{i, D|_{1, \dots, i-1}, D_i}(I)}_{\substack{\text{probability for algorithm }A^{\infty, \epsilon} \\ \text{to select decision }D_i \text{ if the }\\ \text{previous decisions taken are }D|_{1, \dots, i-1}}} f(I, D)$$

    The probabilities taken by the algorithm $A^{\infty, \epsilon}$ are close to the probabilities taken by algorithm $A^j$, thus we write:
    $$ p^{\infty, \epsilon}_{i, D|_{1, \dots, i-1}, D_i}(I) = p^{j}_{i, D|_{1, \dots, i-1}, D_i}(I) + \delta^j_{i, D|_{1, \dots, i-1}, D_i}(I)$$
    with $|\delta^j_{i, D|_{1, \dots, i-1}, D_i}(I)|\leq \Delta_i$

    Then:
    $$ \mathbb E(f(I, A^{\infty, \epsilon}(I))) = \sum_{d^1 \in \mathcal D}(p^{j}_{1, \emptyset, d^1}(I) + \delta^j_{i, \emptyset, d^1}(I))\sum_{d^2 \in \mathcal D}(p^{j}_{2, (d^1), d^2}(I) + \delta^j_{2, (d^1), d^2}(I)) \;\dots\;  f(I, (d^1, d^2, \dots, d^s))$$

    Moreover, without the terms $\delta$, then the expression becomes the expected value of algorithm $A^j$ on $I$:
    $$ \mathbb E(f(I, A^{j}(I))) = \sum_{d^1 \in \mathcal D}p^{j}_{1, \emptyset, d^1}(I)\sum_{d^2 \in \mathcal D}p^{j}_{2, (d^1), d^2}(I)\;\dots\;  f(I, (d^1, d^2, \dots, d^s))$$

    Let $D\in \mathcal D^s$, and define $e^{\infty, \epsilon}_{s, D}(I) = f(I, D)$. Define then inductively for $0\leq i< s$ and for $D\in \mathcal D^{i}$:
    $$e^{\infty, \epsilon}_{i, D}(I) = \sum_{d^{i+1}\in\mathcal D} (p^{j}_{i+1, D, d^{i+1}}(I) + \delta^j_{i+1, D, d^{i+1}}(I))e^{\infty, \epsilon}_{i+1, D\oplus d^{i+1} }(I)$$
    where $\oplus$ in the above equation is the concatenation of sequences. The terms $e^{\infty,\epsilon}$ may be seen as a sub-part of the sum to compute the expected value of algorithm $A^{\infty, \epsilon}$:
    $$ \mathbb E(f(I, A^{\infty, \epsilon}(I))) = \underbrace{\sum_{d^1 \in \mathcal D}(p^{j}_{1, \emptyset, d^1}(I) + \delta^j_{i, \emptyset, d^1}(I))\underbrace{\sum_{d^2 \in \mathcal D}(p^{j}_{2, (d^1), d^2}(I) + \delta^j_{2, (d^1), d^2}(I)) \;\dots\;  \underbrace{f(I, (d^1, d^2, \dots, d^s))}_{e^{\infty, \epsilon}_{s, (d^1, \dots, d^s)}}}_{e^{\infty, \epsilon}_{1, (d^1)}}}_{e^{\infty, \epsilon}_{0, \emptyset}}$$

    We define similarly $e^{j}_{s, D}(I) = f(I, D)$ and inductively for $0\leq i< s$ and for $D\in \mathcal D^{i}$:
    $$e^{j}_{i, D}(I) = \sum_{d^{i+1}\in\mathcal D} p^{j}_{i+1, D, d^{i+1}}(I)\cdot e^{j}_{i+1, D\oplus d^{i+1} }(I)$$
    
    Now, we will prove that:
    $$|\mathbb E(f(I, A^{\infty, \epsilon}(I))) - \mathbb E(f(I, A^{j}(I)))| = |e^{\infty, \epsilon}_{0, \emptyset} - e^{j}_{0, \emptyset}|\leq \epsilon$$

    To do that, we show the following by backwards induction on $i$:
    $$\forall D\in\mathcal D^i,\;\;|e^{\infty, \epsilon}_{i, D} - e^{j}_{i, D}|\leq \frac{\epsilon}{2^i}$$
    Clearly, the statement is true when $i=s$, as $\forall D\in\mathcal D^s,\;\;e^{\infty, \epsilon}_{s, D} = e^{j}_{s, D} = f(I, D)$. Suppose that the statement is true for some $0\leq i\leq s$, let us show that it holds as well for $i-1$.

    $$\forall D\in\mathcal D^{i-1},\;\;|e^{\infty, \epsilon}_{i-1, D} - e^{j}_{i-1, D}| = |\sum_{d^{i}\in\mathcal D} p^j_{i, D, d^i}(e^{\infty, \epsilon}_{i, D\oplus d^{i}} - e^{j}_{i, D\oplus d^i}) + \delta^j_{i, D, d^{i}}\cdot e^{\infty, \epsilon}_{i, D\oplus d^{i}} |$$
    By induction, we know that $|e^{\infty, \epsilon}_{i, D\oplus d^{i}} - e^{j}_{i, D\oplus d^i}|$ is smaller than $\frac{\epsilon}{2^i}$. Moreover, $\sum_{d^{i}\in\mathcal D} p^j_{i, D, d^i} = 1$ since the term $p^j_{i, D, d^i}$ represents a probability. Additionally, $\delta^j_{i, D, d^{i}}\leq \Delta_i = \frac{\epsilon}{2^i m^2}$, and $e^{\infty, \epsilon}_{i, D\oplus d^{i}}$ is a mean over values of $f$ and is thus smaller than $m$, as $|f|\leq m$ (the maximum highest load of a bin is $m$ when all items are put into the same bin).
    Putting everything together:
    $$\forall D\in\mathcal D^{i-1},\;\;|e^{\infty, \epsilon}_{i-1, D} - e^{j}_{i-1, D}| = |\underbrace{\sum_{d^{i}\in\mathcal D} p^j_{i, D, d^i}}_{\text{sums to 1 }}\underbrace{(e^{\infty, \epsilon}_{i, D\oplus d^{i}} - e^{j}_{i, D\oplus d^i})}_{\leq \frac{\epsilon}{2^i}} + \underbrace{\sum_{d^{i}\in\mathcal D}}_{m\text{ terms}} \underbrace{\delta^j_{i, D, d^{i}}}_{\leq \frac{\epsilon}{2^i m^2}}\cdot \underbrace{e^{\infty, \epsilon}_{i, D\oplus d^{i}}}_{\leq m} | \leq \frac{\epsilon}{2^{i-1}}$$

    Hence:
    $$|\mathbb E(f(I, A^{\infty, \epsilon}(I))) - \mathbb E(f(I, A^{j}(I)))| = |e^{\infty, \epsilon}_{0, \emptyset} - e^{j}_{0, \emptyset}|\leq \epsilon$$

    In other words, for any item sequence $I\in\bigcup_{k\in\mathbb N}\mathcal I_{2^k}$ the expected highest load of a bin of algorithm $A^{\infty, \epsilon}$ is at most the expected highest load of a bin of algorithm $A^j$ plus $\epsilon$, with $j\in G(I)$.

    \textbf{3. Extension of the algorithm to process any item sequence:}

    The algorithm $A^{\infty, \epsilon}$ may be extended to be able to process any input item sequence (so far, it has been only defined on item sequences from $\bigcup_{k\in\mathbb N}\mathcal I_{2^k}$). To do so, let $\delta>0$; we construct yet another online algorithm, $A^{\infty, \epsilon, \delta}$ which may process any input sequence and has worst-case expected performance close, up to $\delta$, from the performance of algorithm $A^{\infty, \epsilon}$. The idea is to round items to nearby numbers of the form $\frac{x}{2^k}$ which may be processed by the algorithm $A^{\infty, \epsilon}$. This is done through the following function $T$:
    $$T^\delta \colon \mathcal I \to \bigcup_{k\in\mathbb {N}} \mathcal I_{2^k}$$
    That function is constructed as a sequence of functions (similarly to online algorithms) where the $i$-th function takes as input only the $i$-th item in the sequence:
    $$T^\delta(I) = (T^\delta_1(I_1), T^\delta_2(I_2), \dots)$$
    Let $i_0 = \lceil \log_2(\frac{2}{\delta}) \rceil$.
    $$T^\delta_i(y) = \frac{\lfloor 2^{i_0 + i} y\rfloor}{2^{i_0+i}}$$ 
    This construction implies:
    $$\forall i\in\mathbb N^*, \; y - T^\delta_i(y) \leq \frac{\delta}{2^i} $$
    Hence, $T^\delta$ verifies:
    $$||I - T^\delta (I)||_{1} \leq \delta$$

    Then, define $A^{\infty, \epsilon, \delta}(I) = A^{\infty, \epsilon}(T^\delta(I))$: in other words, any received item if first "rounded" down to some item of the form $\frac{x}{2^k}$ for some $x$ and $k$ integers, then the algorithm $A^{\infty, \epsilon}$ is applied. Doing so induces an error of at most $\delta$ (since items are considered smaller than what they really are by the algorithm $A^{\infty, \epsilon}$). Since items are always rounded down, the sequence $T(I)$ is still a valid sequence since it does not violate the constraint of the online bin stretching problem that its offline optimum must be lesser than 1.

    The extended algorithm $A^{\infty, \epsilon, \delta}$ is still comprised of measurable transformations, as $(T^\delta_i)^{-1}$ is by construction a finite union of  segments. Hence, we have constructed a valid online randomized algorithm, with performance of at most $\epsilon + \lim\limits_{k\to +\infty} v^*_{2^k} + \delta$. The performance of any algorithm is an upper bound on the optimal value of the game. Hence: $$\forall \epsilon>0,\;\forall \delta>0,\;\lim\limits_{k\to +\infty} v^*_{2^k} \leq v^* \leq \lim\limits_{k\to +\infty} v^*_{2^k} + \epsilon + \delta $$
    This concludes the proof that the linear programming method converges towards the optimal competitive ratio for the online bin stretching problem. 

\end{proof}

\bibliography{refs}

\end{document}